\documentclass[11pt,reqno]{amsart}

\theoremstyle{plain}
\newtheorem{theorem}{Theorem}

\newtheorem{lemma}[theorem]{Lemma}
\newtheorem{corollary}[theorem]{Corollary}
\newtheorem{fact}[theorem]{Fact} 
\theoremstyle{definition}

\theoremstyle{remark}
\newtheorem{remark}[theorem]{Remark}
\newtheorem*{Acknowledgement}{Acknowledgement}
\catcode`@=11 \@mparswitchfalse  

\newcommand{\hsp}[1]{\hskip #1 pt}
\newcommand{\vsp}[1]{\vskip #1 pt}
\newcommand{\e}{\varepsilon}
\newcommand{\ai}{asymptotically isometric}
\newcommand{\tr}{\mathcal T}
\newcommand{\pe}{\widehat{\imath}\:}
\newcommand{\ft}{\tilde{f}}
\newcommand{\at}{\widetilde{a}}
\newcommand{\gn}{\gamma_n}
\newcommand{\gnn}{\e}  
\newcommand{\defeq}{\hsp{5} \overset{\text{\rm def}}{=} \hsp{5}}
\newcommand{\arrowi}{\overset{\text{isomorphic}}{\vector(1,0){60}}}
\newcommand{\arrowali}{\overset{\text{almost isometric}}{\vector(1,0){60}}}
\newcommand{\arrowasi}{\overset{\text{asymp.\ isometric}}{\vector(1,0){60}}}
\newcommand{\arrowt}{\overset{\text{isometric}}{\vector(1,0){60}}}

\begin{document}

\date{July 4, 1999}

\title{Dual Banach spaces which contain an isometric copy of $L_1$}
\author{S.\ J.\ Dilworth, Maria Girardi and J.\ Hagler}
\address{Department of Mathematics, University of South Carolina, 
         Columbia, SC 29208, U.S.A.} 
\email{dilworth@math.sc.edu and girardi@math.sc.edu}
\address{Department of Mathematics and Computer Science,
         University of Denver, Denver, CO 80208, U.S.A.}
\email{jhagler@cs.du.edu}
\subjclass{Primary 46B04. Secondary 46B20.}
\thanks{Maria Girardi is supported by NSF grant DMS--9622841}
\keywords{Banach spaces, \ai\ copies of $\ell_1$, 
isometric copy of $L_1$}

\begin{abstract} 
A Banach space contains 
\ai\ copies of~$\ell_1$ 
if and only if 
its dual space contains 
an isometric copy of $L_1$.
\end{abstract}

\maketitle

\baselineskip 12.7 pt 
 
\section{Introduction}\label{S:intro}

The duality between a Banach space containing 
a  `nice' copy of $\ell_1$ 
and its dual space containing a  `nice'
copy of $L_1$ is summarized 
in the diagram below. 
Each upward implication follows straight from the 
definitions and  
the absence of a downward arrow indicates that 
the corresponding implication does not hold in general.
 
$$
\begin{picture}(300,250)
\put(0,0){\makebox(0,0){$\ell_1$}}
\put(50,8){\makebox(0,0){$\arrowt$}}
\put(100,0){\makebox(0,0){$X$}}
\put(50,44){\makebox(0,0){\vector(0,1){20}}}
\put(0,80){\makebox(0,0){$\ell_1$}}
\put(50,88){\makebox(0,0){$\arrowasi $}}
\put(100,80){\makebox(0,0){$X$}}
\put(50,124){\makebox(0,0){\vector(0,1){20}}}
\put(0,160){\makebox(0,0){$\ell_1$}}
\put(50,168){\makebox(0,0){$\arrowali $}}
\put(100,160){\makebox(0,0){$X$}}
\put(47,210){\makebox(0,0){\vector(0,1){20}}}
\put(53,210){\makebox(0,0){\vector(0,-1){20}}}
\put(0,240){\makebox(0,0){$\ell_1$}}
\put(50,248){\makebox(0,0){$\arrowi $}}
\put(100,240){\makebox(0,0){$X$}}

\put(200,0){\makebox(0,0){$L_1$}}
\put(250,8){\makebox(0,0){$\arrowt$}}
\put(300,0){\makebox(0,0){$X^*$}}
\put(250,80){\makebox(0,0){\vector(0,1){40}}}
\put(200,160){\makebox(0,0){$L_1$}}
\put(250,168){\makebox(0,0){$\arrowali $}}
\put(300,160){\makebox(0,0){$X^*$}}
\put(247,210){\makebox(0,0){$\vector(0,1){20}$}}
\put(253,210){\makebox(0,0){$\vector(0,-1){20}$}}
\put(200,240){\makebox(0,0){$L_1$}}
\put(250,248){\makebox(0,0){$\arrowi $}}
\put(300,240){\makebox(0,0){$X^*$}}

\put(150,242){\makebox(0,0){\vector(1,0){30}}}
\put(150,242){\makebox(0,0){\vector(-1,0){30}}}
\put(150,40){\makebox(0,0){\vector(1,-1){40}}}
\put(150,40){\makebox(0,0){\vector(-1,1){40}}}
\end{picture}
$$
\vsp{5} 
 
The investigation of this duality began when 
Pe\l czy\'nski~\cite{P} proved that if  $X$ contains
an isomorphic copy of  $\ell_1$ then $X^*$ contains
an  isomorphic copy of~$L_1$. He also proved the
converse result under a technical assumption which was later
removed by Hagler~\cite{H}. 
Earlier, James \cite{J} had shown that if~$X$ 
contains $\ell_1$ isomorphically     
then $X$ contains $\ell_1$ almost 
isometrically. 
Recently, Dowling, N.\ Randrianantoanina and
Turett~\cite{DRT} proved that a dual Banach space contains
almost isometric copies of $L_1$ whenever it contains
isomorphic copies of $L_1$ 
(see also \cite[Corollary~2.32]{H1} for this result).
The main result of this paper, Theorem~\ref{th: mainresult}, 
shows that $X$ contains \ai\ copies of $\ell_1$ if 
and only if $X^*$ contains $L_1$ isometrically. 
In the real case,  this is a hitherto unpublished result of Hagler \cite[Theorem~2.2]{H1}.

\section{Notation and Terminology}\label{S:notation}

Henceforth, all Banach spaces are 
either real or complex. 
$X$, $Y$, and~$Z$ will denote   
arbitrary (infinite-dimensional) Banach spaces.    
Let 
$C(K)$  be the space of continuous functions
on  some compact Hausdorff space $K$, 
let~ 
$L_1$ be the space of Lebesgue--integrable functions on [0,1],
and let 
$\ell_p(\Gamma)$ be the space of scalar-valued functions on the set $\Gamma$ 
with finite $\| \cdot \|_p$-norm where~$1\le p \le \infty$,    
all with their usual norms. 
Let $\Delta$ be the Cantor set,~$\ell_p$ be~$\ell_p(\mathbb N)$, and~
$C$ be $C([0,1])$.

The  concept   
of {\it \ai\  copies of $\ell_1$\/} was introduced 
by Hagler \cite[pg.~14]{H1}. 
It  was    
revitalized recently by   Dowling and Lennard    
in fixed point theory \cite{DL}.
A Banach space contains {\it \ai\ copies of $\ell_1$} 
provided it satisfies one
(hence all)  of the conditions 
in the lemma below. 
 
\begin{lemma}\label{ai}
For a Banach space $X$, the following are equivalent.
\begin{enumerate}
\item[\rm{(A1)}]
There exist  
a null sequence $(\varepsilon_n)$ of positive  numbers 
less than one  and   
a sequence $(x_n)$ in $X$ such that
\begin{equation*}
\sum_{n=1}^m (1-\varepsilon_n)|a_n| \le \left\|\sum_{n=1}^m
a_n x_n\right\| \le \sum_{n=1}^m |a_n|  
\end{equation*}
for each finite sequence $(a_n)_{n=1}^m$ of scalars.
\item[\rm{(A2)}]
There  exist  a null sequence $(\varepsilon_n)$ of positive  numbers
and a sequence $(x_n)$ in $X$ such that
\begin{equation*} 
\sum_{n=1}^m |a_n| \le \left\|\sum_{n=1}^m
a_n x_n\right\| \le \sum_{n=1}^m (1+\varepsilon_n) |a_n|  
\end{equation*} 
for each finite sequence $(a_n)_{n=1}^m$ of scalars.
\item[\rm{(A3)}] 
There exist  
a null sequence $(\varepsilon_n)$ of positive  numbers
and a sequence $(x_n)$ in $X$ such that
\begin{equation*} 
\sum_{n=k}^m |a_n| \le \left\|\sum_{n=k}^m
a_n x_n\right\| \le (1+\varepsilon_k) \sum_{n=k}^m  |a_n| 
\end{equation*}    
for each finite sequence $(a_n)_{n=k}^m$ of scalars 
and $k\in \mathbb N$.
\end{enumerate}
\end{lemma}
\noindent 
The proof of this lemma is elementary  
(cf.~\cite[Theorem~1.7]{DLT} for further equivalent formulations). 
Note that   each condition is equivalent to 
the variant obtained by replacing 
{\it  `There exist a'\/}  
 by {\it  `For  each'\/} 
and  {\it `and'\/} by {\it `there exists\rm{.}'\/}
A sequence $(x_n)$ satisfying one of the  conditions in the lemma 
is called  {\it an \ai\ copy of ~$\ell_1$\/}.   
See ~\cite{DLT}   for a splendid survey of this topic and 
its applications to fixed point theory.  
 
The proof of James's theorem \cite{J} for $\ell_1$ shows that 
if $X$ contains  $\ell_1$ almost isometrically, 
then for each    
null sequence $(\varepsilon_n)$ of positive  numbers
there exists a sequence $(x_n)$ in $X$ such that
\begin{equation*}
(1-\varepsilon_k)\sum_{n=k}^m  |a_n| \le \left\|\sum_{n=k}^m
a_n x_n\right\| \le \sum_{n=k}^m |a_n|  
\end{equation*}
for  each finite sequence $(a_n)_{n=k}^m$ of scalars and $k\in\mathbb N$.
Indeed, the line 
between containing $\ell_1$ 
almost isometrically  
and asymptotically isometrically   is very fine.

A sequence $(x_n)$ in a Banach space $X$ is 
a ($1+\e$){\it-perturbation of an isometric copy of $\ell_1$\/} 
(for short, a ($1+\e$){\it-p.i.\ $\ell_1$-sequence\/})   
provided that there 
exist a Banach space $Y$,
a  linear isometric embedding 
$T \colon X \to Y$,  
and  a sequence
$(y_n)$ in $Y$ such that~ $(y_n)$ is 
isometrically equivalent to the unit  vector basis 
of $\ell_1$ and 
$
\left\|  y_n - Tx_n \right\| \leq \e  
$ 
for each $n\in\mathbb N$. 
If furthermore 
\begin{equation*}
\left\|  y_n - Tx_n \right\| 
\defeq \e_n \hsp{4} 
\xrightarrow{n\to\infty} \hsp{4} 0 \ ,  
\end{equation*}
then $(x_n)$ is a {\it perturbation of an isometric copy of $\ell_1$\/} 
(for short, a {\it p.i.\ $\ell_1$-sequence\/}) 
with respect to $(\e_n)$. 
Note that if $X$ is separable  then $Y$ may 
be taken to be separable.

If $X$ is a Banach space, then 
$X^*$ is its dual space,  
$B(X)$ is its closed unit ball, and 
$S(X)$ is its unit sphere. 
The closed linear span of a subset~$A$ of~$X$ is~$[A]$.
If $Y$ is a subspace of $X$ then $\pi\colon X \to X/Y$ 
is the natural quotient mapping.  

For a    surjective bounded linear operator $T \colon X \to Z$, 
the corresponding bounded linear operator $T_q$ is defined by 
the following (commutative) diagram.   
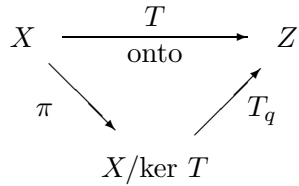
\vsp{3}
$$
\begin{picture}(100,60)
\put(0,50){\makebox(0,0){$X$}}
\put(100,50){\makebox(0,0){$Z$}}
\put(50,0){\makebox(0,0){$X/{\text{ker}~ T}$}}
\put(15,50){\vector(1,0){70}}
\put(50,58){\makebox(0,0){$T$}}
\put(50,44){\makebox(0,0){\text{onto}}}
\put(10,40){\vector(1,-1){25}}
\put(5,20){$\pi$}
\put(65,13){\vector(1,1){25}}
\put(85,20){$T_q$}
\end{picture}
$$ 
\vsp{10}\noindent
The operator $T$ is called an {\it isometric quotient mapping\/}   
provided $T_q$ is an isometry, which is 
the case if and only if $T^*$ is an isometric embedding.  
If~$S \colon X \to Z$ is an 
isomorphic embedding,  
then $S^*$ is an isometric quotient mapping if and only if  
$S$ is an isometric embedding.

All notation and terminology, not otherwise explained, 
are as in~\cite{LT}.

\section{Main Result}\label{S:mainresult}

Theorem~\ref{th: mainresult}, the main result of this paper,   
may be viewed as the
isometric version of the theorems of
Pe\l czy\'nski and Hagler.

\begin{theorem}\label{th: mainresult} 
For a Banach space $X$, the 
following are equivalent. 
\begin{enumerate}
\item[\rm{(a)}] 
$X$ contains \ai\ copies of $\ell_1$.
\item[\rm{(b)}] 
$X$ contains a perturbation of an isometric copy of $\ell_1$.
\item[\rm{(c)}]
$\ell_1$ is linearly isometric to a quotient space of a subspace $X$.  
\item[\rm{(d)}]
$L_1$ is linearly isometric to a subspace of $X^*$.
\item[\rm{(e)}]
$C^*$ is linearly isometric to a subspace of $X^*$.
\item[\rm{(f)}]
$X^*$ contains an infinite set $\Gamma$  which
is isometrically equivalent to the usual basis of $\ell_1(\Gamma)$
and which is dense-in-itself in the weak-star topology on $X^*$.  
\end{enumerate}
And if $X$ is separable, then the following is   equivalent 
to each of the above conditions. 
\begin{itemize}
\item[\rm{(g)}] 
$C(\Delta)$ is isometric to a quotient space of $X$. 
\end{itemize}
\end{theorem}

\noindent
Recall that a subset  $K$  of a topological space is 
dense-in-itself 
if $K$ has no isolated points in the  relative topology. 
Our proof of 
Theorem~\ref{th: mainresult}  
uses the following results.   

\begin{lemma} 
\label{lem:pi}
If $(x_n)$ is a p.i.\ $\ell_1$-sequence, 
then $(\lambda_n x_n)$ is an \ai\ copies of $\ell_1$  
satisfying~(A2) 
for some suitable choice of scalars~$(\lambda_n)$. 
Conversely, an \ai\ copies of $\ell_1$ satisfying~(A2) 
is a p.i.~$\ell_1$-sequence.    
\end{lemma} 

\begin{proof}
Let  $(\tilde{x}_n)$  
be a p.i.\  $\ell_1$-sequence  
with respect to $(\tilde{\e}_n)$. Then     
$$
\sum_{n=1}^m (1-\tilde{\e}_n)|a_n| \le \left\|\sum_{n=1}^m
a_n \tilde{x}_n\right\| \le \sum_{n=1}^m (1+\tilde{\e}_n)|a_n|  \ . 
$$ 
for each finite sequence $(a_n)_{n=1}^m$ of scalars.  
Define  
$$
\e_n \defeq \frac{1+\tilde{\e}_n}{1-\tilde{\e}_n} ~-~1
\hsp{30}\text{and}\hsp{30}
x_n \defeq \frac{\tilde{x}_n}{1-\tilde{\e}_n} \hsp{10} . 
$$ 
Then $(\e_n)$ and $(x_n)$ satisfy (A2); 
thus, $(x_n)$ is an \ai\ copy of $\ell_1$.  
 
Conversely,  
let  $(\e_n)$ and $(x_n)$ 
satisfy~(A2).  Then $(x_n)$ is a p.i.\ sequence.       
To see this,  let~$X_0~=~[ x_n ]$
and   
$$
W = \left\{ \left( w_n \right)_{n=1}^\infty \colon 
w_n \in \mathbb C 
\text{\hsp{6}and\hsp{6}} |w_n|=1 
\text{\hsp{6}for each\hsp{6}} n\in \mathbb N \right\}. 
$$ 
For each $\omega = (w_n) \in W$, define $f_\omega \in B(X_0^*)$ 
by $f_\omega(x_n) = w_n $;  for indeed, 
\begin{equation*}
\left| f_\omega \left(\sum_{n=1}^m a_n x_n \right) \right| 
\hsp{4}=\hsp{4} 
\left| \sum_{n=1}^m a_n w_n  \right| 
\hsp{4}\le\hsp{4} 
\sum_{n=1}^m \left| a_n \right|
\hsp{4}\le\hsp{4} 
\left\| \sum_{n=1}^m a_n x_n \right\|_{X_0}  
\end{equation*}
for each finite sequence $(a_n)_{n=1}^m$ of scalars.   
For each $\omega\in W$, let $\ft_\omega \in B(X^*)$ 
be a norm-preserving Hahn-Banach extension of $f_\omega$. 

Let  
\begin{equation*}
Y \defeq C\left( B\left(X^*\right) , \sigma\left(X^*, X\right)\right) \ ,  
\end{equation*} 
endowed with the usual sup norm,  
and consider the isometric embedding 
\begin{equation*}
 T \hsp{2}\colon X  \hsp{4} \to\hsp{4}  Y 
\end{equation*} 
given by
\begin{equation*}
 (Tx)(x^*) \defeq x^*(x) \ .  
\end{equation*}   
Let $y_n \in B(Y)$ be the  
`truncation' of $T x_n$; specifically,
\begin{equation}
\label{eq: trunc} 
y_n (x^*) \hsp{5}=\hsp{5}
\begin{cases}
\hsp{5} ( T x_n )(x^*) &\hsp{8}\text{if}\hsp{12} 
|(T x_n )(x^*) | \le 1 \\  
\hsp{5} \dfrac{(T x_n )(x^*)}{|(T x_n )(x^*)|} 
&\hsp{8}\text{if}\hsp{12} 
|(T x_n )(x^*) |  > 1 \ . 
\end{cases}
\end{equation} 

For each $n\in \mathbb N$,      
condition~(A2) gives that 
$  
\left\|   x_n \right\| \le 1 + \e_n, 
$  
and so by~\eqref{eq: trunc}  
\begin{equation*}
\| y_n - Tx_n \|_Y 
~\le~ \e_n  \ .  
\end{equation*}
Since for each $n\in\mathbb N$ and $\omega = (w_j) \in W$ 
\begin{equation*}
(Tx_n)(\ft_\omega) 
\hsp{3}=\hsp{3} 
\ft_\omega (x_n)
\hsp{3}=\hsp{3} 
 f_\omega (x_n) 
\hsp{3}=\hsp{3}
w_n  
\hsp{3}=\hsp{3}
y_n(\ft_\omega)   \ ,   
\end{equation*}
it follows that  
\begin{equation*}
\left\| \sum_{n=1}^m a_n y_n \right\|_Y 
\hsp{3}\geq\hsp{3}
\sup_{\omega\in W} 
\left| \sum_{n=1}^m a_n y_n(\ft_\omega) \right|
\hsp{3}=\hsp{3}
\sup_{ (w_n)\in W} 
\left| \sum_{n=1}^m a_n  w_n \right|
\hsp{3}=\hsp{3}
\sum_{n=1}^m \left| a_n \right|   \ .  
\end{equation*}
for  each   finite sequence $(a_n)_{n=1}^m$ of 
scalars.  
Also, $\left\| y_n \right\| \le 1$ for each $n\in\mathbb N$.    
Thus  $(y_n)$ is isometrically equivalent 
to the unit vector basis of~$\ell_1$.
\end{proof}

\begin{remark}
Minor  modifications to the above proof give    
an isomorphic version of  Lemma~\ref{lem:pi}.  
Indeed,  if  $(\tilde{x}_n)$   
be a ($1+\tilde{\e}$)-p.i.\  $\ell_1$-sequence   
with $\tilde{\e}<1$  
and 
$$
\e \defeq \frac{1+\tilde{\e}}{1-\tilde{\e}} ~-~1
\hsp{30}\text{and}\hsp{30}
x_n \defeq \frac{\tilde{x}_n}{1-\tilde{\e}} \hsp{10} ,  
$$ 
then 
\begin{equation}
\label{eq:ipi}
\sum_{n=1}^m|a_n| \le \left\|\sum_{n=1}^m
a_n x_n\right\| \le (1 + \e) \sum_{n=1}^m |a_n|  
\end{equation}
for each finite sequence $(a_n)_{n=1}^m$ of scalars. 
Conversely,  
if $( x_n )$ satisfies~\eqref{eq:ipi} 
for each finite sequence $(a_n)_{n=1}^m$ of scalars,  
then $(x_n)$ is a  ($1+\e$)-p.i.\  $\ell_1$-sequence.    
\end{remark}

\begin{lemma} 
\label{lem: sep} 
If $X$ satisfies (f) of Theorem~\ref{th: mainresult}, 
then there  exists a separable subspace~$X_0$ of $X$  
and a countable subset $\Gamma^\prime$ of $X_0^*$ which  
satisfies~(f) of Theorem~\ref{th: mainresult}. 
\end{lemma}

\begin{proof}
Let $X$ be a Banach space satisfying (f) of Theorem~\ref{th: mainresult}. 
We shall inductively construct
a sequence~$(\Lambda_i)$ of  countably infinite subsets of $\Gamma$
and  
a sequence $(Z_i)$ of separable subspaces of~$X$   
which satisfy, for each~$n\in \mathbb N$,  
\begin{enumerate}
\item[\rm{(1)}] 
$Z_n \subset Z_{n+1}$, 
\item[\rm{(2)}]
 $Z_n$ norms $[\:\cup_{i=1}^n \Lambda_i\:]$, 
\hsp{10} i.e., if $x^* \in \left[\:\cup_{i=1}^n \Lambda_i\:\right]$   
then 
\begin{equation*}
\left\| x^* \right\| = \sup_{z \in B(Z_n)}|x^*(z)| \ , 
\end{equation*}
\item[\rm{(3)}]
$\cup_{i=1}^n \Lambda_i$ \, is contained in 
the $Z_n$-cluster points of $\Lambda_{n+1}$, 
\\  i.e., 
if $x^* \in \cup_{i=1}^n \Lambda_i$  \, 
and $(w_i)_{i=1}^k$ are from $Z_n$ and $\e>0$   
then 
\begin{equation*}
\left\{ y^* \in \Lambda_{n+1} \colon \left| \left( y^* - x^* \right) 
\left( w_i \right) \right| < \e \hsp{3}\text{for}\hsp{3} 
1 \le i \le k \right\} \setminus \left\{ x^* \right\} 
\not = \emptyset \ . 
\end{equation*} 
\end{enumerate} 
For the first step of the induction choose a countably infinite subset
$\Lambda_1$ of~$\Gamma$ and find 
a separable subspace $Z_1$ of $X$ which satisfies (2).  
Suppose that we have chosen $(\Lambda_i)_{i=1}^n$ and $(Z_i)_{i=1}^n$ 
satisfying the three conditions.  
Since~$Z_n$  is separable and 
elements of $\Gamma$ are of norm one,  
there is a countable
subset~$\Lambda_{n+1}$ of $\Gamma$ 
satisfying~(3).   
Next we find a separable
subspace $Z_{n+1}$ of $X$ which 
satisfies~(1) and ~(2).  
This completes the inductive step.

Now let $X_0=[\cup_{n=1}^\infty Z_n]$ and 
\begin{equation*}
\Gamma^\prime \defeq
\left\{ x^*\vert_{X_0} \colon x^* \in  \cup_{n=1}^\infty \Lambda_n \right\} \ . 
\end{equation*} 
Condition~(2) gives that $\Gamma^\prime$ is 
isometrically equivalent to the
usual basis of~$\ell_1(\Gamma^\prime)$.  
Conditions~(1) and~(3), along with the fact that 
$\Gamma \subset S(X^*)$,  
give  that~$\Gamma^\prime$ is
dense-in-itself in the weak-star topology on $X_0^*$.
\end{proof}

\begin{fact} 
{\rm (cf.~{\cite[Lemma~4]{HS}})} 
\label{topology} 
Let $N$ and $M$ be compact Hausdorff spaces
with~$M$ perfect and suppose that $\phi \colon N \rightarrow M$ is continuous
and onto. Then there exists a subset~$Q$ of $N$ such 
that $Q$ is dense-in-itself and
$\phi|_Q \colon Q \to M$ is a bijection. 
\end{fact}

\begin{fact}
{\rm (Haskell P.~Rosenthal {\cite[Proposition~3 and its Remark~2]{R}})}
\label{thm: hpr} 
Let $X_0$ be a  
separable Banach space satisfying (f) of 
Theorem~\ref{th: mainresult}. Then  there exists 
$$
K \subset B(X^*_0)  \ , 
$$ 
which is homeomorphic to $\Delta$, 
such that  the restriction operator    
\vsp{8}\noindent
$$
R \hsp{2}\colon X_0 \hsp{4} \to\hsp{4} C(K) 
$$ 
given by 
$
 (Rx_0)(x_0^*) 
\hsp{3}=\hsp{3} 
x_0^*(x_0)
$
is an isometric quotient mapping. 
\end{fact}

\begin{proof}[Proof of Theorem~\ref{th: mainresult}]
We shall assume that $X$ is a complex Banach space as the proof in
the real case is easier.

The equivalence of (a) and (b) follows directly 
from Lemma~\ref{lem:pi}.

To see that (a) implies (c), 
let $(\e_n)$ and  $(x_n)$ be   sequences   satisfying condition~(A3).  
Partition $\mathbb N$ into infinite sets $\{ J_n \}_{n\in \mathbb N}$  
and let   
$ 
T \colon [ x_n ]  \to \ell_1 
$  
be the bounded linear operator that maps 
$x_j$ to the  $n^{\text{th\/}}$ unit vector of $\ell_1$ 
when $j\in J_n$. 
Then    $T$ is an 
isometric quotient mapping.  

To see that (c) implies (a), 
let 
$$
T \colon {X_0}/{X_1} \to \ell_1 
$$ 
be an isometry from a quotient space of a subspace $X_0$ of $X$ 
onto $\ell_1$.  Fix a 
null sequence $(\e_n)$ 
of positive 
numbers.  Find a sequence $(x_n)$  in $X_0$  
such that~$T(x_n + X_1)$ is the  $n^{\text{th\/}}$ unit vector of $\ell_1$
and
$$ 
1 \le \left\| x_n \right\|_X \le 1+\e_n  \ . 
$$ 
Then $(\e_n)$ and  $(x_n)$ satisfy~(A2).  

To see that (a) implies (e), 
let $(\e_n)$ and  $(x_n)$ satisfy~(A2).  
We shall define the  bounded linear operators
in the (commutative) diagram below 
\vsp{5}\noindent 
\begin{equation} 
\begin{picture}(230,50)
\put(0,0){\makebox(0,0){$[x_n] \defeq $}}
\put(30,0){\makebox(0,0){$Y$}}
\put(80,4){\makebox(0,0){\vector(1,0){60}}}
\put(80,6){\makebox(0,0){{$T$}}}
\put(130,0){\makebox(0,0){{$C$}}}
\put(180,4){\makebox(0,0){{\vector(1,0){60}}}}
\put(180,6){\makebox(0,0){{$\pe$}}}
\put(230,0){\makebox(0,0){{$C^{**}$}}}
\put(30,50){\makebox(0,0){{$X$}}}
\put(30,25){\makebox(0,0){{\vector(0,1){30}}}}
\put(20,25){\makebox(0,0){{$j$}}} 
\put(130,30){\makebox(0,0){{\vector(4,-1){160}}}}
\put(130,42){\makebox(0,0){{$\widetilde T$}}}
\end{picture}
\label{eq: cd}
\end{equation}
\vsp{5}\noindent 
as follows.  
Let $(z_n)$ be dense in the unit sphere of $C$. 
Define $T$ by $Tx_n~=~z_n$.  
Condition~(A2) gives that  
$T$ is a surjective norm-one  bounded 
linear operator. 
Furthermore, $T$ is an    
isometric quotient mapping for 
if~$f\in S(C)$  then there is a subsequence $(z_{k_n})$ 
converging in norm to $f$ and so  
$$
1~\le~\left\| T^{-1}_q f \right\|_{Y/{\text{ker~}T}}~\le~
\varliminf_{n\to\infty} \left\| x_{k_n} \right\|_{X} 
~ = ~1 \ . 
$$
Let~$j$ be the natural embedding   
and let $\pe$  be  
 the canonical isometric embedding  
given by point evaluation.  
Since~$C^{**}$ has the  Hahn-Banach Extension Property, 
$T$ admits a norm-preserving extension~$\widetilde T$. 
Dualizing gives the commutative diagram 
\vsp{5} 
$$
\begin{picture}(260,60)
\put(0,50){\makebox(0,0){$C^*$}}
\put(21,54){\makebox(0,0){{\vector(1,0){20}}}}
\put(21,60){\makebox(0,0){{$h$}}}
\put(50,50){\makebox(0,0){$C^{***}$}}
\put(100,54){\makebox(0,0){{\vector(1,0){60}}}}
\put(100,60){\makebox(0,0){{$\pe^*$}}}
\put(150,50){\makebox(0,0){$C^*$}}
\put(200,54){\makebox(0,0){{\vector(1,0){60}}}}
\put(200,60){\makebox(0,0){{$T^*$}}}
\put(250,50){\makebox(0,0){$Y^*$}} 
\put(250,0){\makebox(0,0){$X^*$}}
\put(150,20){\makebox(0,0){{\vector(4,-1){160}}}}
\put(180,25){\makebox(0,0){{$\widetilde{T}^*$}}}
\put(250,25){\makebox(0,0){{\vector(0,1){20}}}}
\put(260,25){\makebox(0,0){{$j^*$}}} 
\end{picture}
$$ 
\vsp{5}\noindent 
where $h$ is the canonical isometric embedding  
given by point evaluation. 
To see that $R \defeq \widetilde{T}^* \: h$   
is the desired isometric embedding, 
let~$\mu\in C^*$.   
Then, since $T^*$ is an isometric embedding and 
$\pe^* \: h$ is the identity mapping, 
$$
\left\| 
\mu \right\|_{C^*}
~=~
\left\| T^*   
\mu  \right\|_{Y^*}
~=~
\left\| j^*  R 
\mu \right\|_{Y^*}
~\le~
\left\|   R 
\mu \right\|_{X^*}
~\le~
\left\| 
\mu \right\|_{C^*} \ . 
$$ 
 
Clearly, (e) implies (d).  

To see that (d) implies (a), let   
$T\colon L_1 \rightarrow X^*$ be an isometric embedding 
and let   
$(\varepsilon_n)$ be a null sequence of positive numbers.
Then $T^* \colon X^{**}\rightarrow L_\infty$ is a weak-star continuous
isometric quotient  mapping.  
By Goldstine's Theorem,   
$$
W \defeq T^*(B(X))
$$ is weak-star dense in~$B(L_\infty)$.
For each $n  \in \mathbb N$, let
$$
F_n  \defeq
\{z^n_j \colon 1 \le j \le M(n)\}
$$
be an $(\varepsilon_n/2)$--net for $\{z\in \mathbb{C}\colon |z|=1\}$.

Let $\tr$ be the tree 
\[
\tr = \mathop{\bigcup}_{n\in\mathbb N} \hsp{3} \tr_n 
\]
where $\tr_n$, the   $n^{\text{th\/}}$-level of  $\tr$, is  
\begin{multline}
\tr_n \defeq \{ (m_0, m_1, m_2, \hsp{3} \ldots \hsp{3} m_{n-1}) 
\in \mathbb N^{\,n} \colon  \\
m_0 = 1
\hsp{5}\text{and}\hsp{5}
1 \le m_j \le M(j) 
\hsp{5}\text{for each}\hsp{5}
j \in \mathbb N \}  \ . 
\notag
\end{multline}
If $\alpha = (m_0, m_1, m_2, \hsp{3} \ldots \hsp{3} m_{n-1}) \in \tr$ 
then   
\[
(\alpha,  j) \defeq 
(m_0, m_1, m_2, \hsp{3} \ldots \hsp{3} m_{n-1},  j)   \ ; 
\]
thus, for each $n\in\mathbb N$  
\[
\tr_{n+1} = \left\{ (a,j) \colon \alpha \in \tr_n 
\hsp{5}\text{and}\hsp{5} 
1 \le j \le M(n) \right\} \ . 
\] 
We will  define inductively, for each $n\in \mathbb N$, a    
collection   $\{ A_\alpha \}_{\alpha \in  \tr_n}$ of disjoint  
sets of positive (Lebesgue) measure  
and a function $f_n\in W$ such that, for each 
$n\in\mathbb N$ and $\alpha\in\tr_n$, 
\begin{equation}
\label{eq: ank1}
\mathop{ \bigcup }_{j=1}^{M(n)} A_{(\alpha,j)} 
\hsp{5}\subset\hsp{5} 
A_\alpha 
\hsp{5}\subset\hsp{5} [0,1] 
\end{equation}
and, for each $1 \le j \le M(n)$,  
\begin{equation}
\label{eq: ank2}
 | f_n   - z^n_j  | < \frac{\e_n}{2} 
\hsp{8}\text{on}\hsp{5}   A_{(\alpha,j)}  \ . 
\end{equation}
To start the induction, 
let $$A_{(m_0)} = [0,1] \ . $$  
For the inductive step, let $n \in \mathbb N$ and 
suppose that we have constructed 
disjoint sets   
\[
\{ A_{\alpha} \colon \alpha \in  \tr_{n} \}  
\]
of  positive  measure.   
For each $\alpha \in  \tr_{n}$,
partition $A_\alpha$ into  sets $\{ D_{(\alpha, j)} \}_{j=1}^{M(n)}$ 
of positive measure.  
Consider the function $g_n \in B(L_\infty)$ defined by 
\begin{equation}
g_n(t) = 
\begin{cases} 
z^n_j & \text{if}\hsp{4} t \in D_{(\alpha, j)}
\hsp{3}\text{and}\hsp{3} \alpha\in\tr_n \\
0 & \text{otherwise} \hsp{5} . 
\end{cases}
\notag
\end{equation} 
Since $W$ is weak-star dense in $B(L_\infty)$ there
exists $f_n \in W$ approximating~$g_n$ closely enough to ensure that the sets
$$
A_{(\alpha, j)} \defeq
\{|f_n- z^n_j| < \varepsilon_n/2\}\cap D_{(\alpha, j)}
$$
all have positive measure. 
This completes the proof of the inductive step. 

For each $n\in \mathbb N$, select
$x_n \in B(X)$ such that $T^*(x_n)=f_n$.
To see that~$(x_n)$ is an \ai\ copy of $\ell_1$, 
let    $(a_n)_{n=1}^m$ be a
finite complex sequence.   
Define $(\at_n)_{n=1}^m$ from   
$\{z\in \mathbb{C}\colon |z|=1\}$  by 
\begin{equation*}
\at_n \hsp{5} = 
\begin{cases} 
\dfrac{\Bar{a}_n}{\left| a_n \right|} 
&\hsp{5}\text{if}\hsp{3} a_n \not = 0 \\
\hsp{6} 1 &\hsp{5}\text{if}\hsp{3} a_n   = 0 \ ;
\end{cases}
\end{equation*}
thus, $a_n  \at_n  = \left| a_n \right| $.  
For each $1 \le n \leq m$, 
find $1\leq j_n \leq M(n)$  so that 
\begin{equation*}
| \at_n -  z^n_{j_n} | \hsp{3} < \hsp{3} \frac{\e_n}{2} \ . 
\end{equation*} 
Then $\alpha \defeq (1, j_1, \dots , j_m) \in \tr_{m+1}$ and so 
by~ \eqref{eq: ank1} and~\eqref{eq: ank2}   
\begin{equation*}
| f_n - z^n_{j_n} | \hsp{3} < \hsp{3} \frac{\e_n}{2} 
\hsp{8}\text{on}\hsp{5}   A_\alpha
\end{equation*} 
for each $1\leq n \leq m$.  Thus  
\begin{align*}
\left\|\sum_{n=1}^m a_n x_n \right\|
\hsp{5}&\ge\hsp{5}
\left\|T^*\left(\sum_{n=1}^m a_n x_n\right)\right\|_{L_\infty}
\hsp{5}=\hsp{5} 
\left\|\sum_{n=1}^m a_n f_n\right\|_{L_\infty}
\\ &\ge\hsp{5}
\left\|\sum_{n=1}^m a_n 
      \at_n 1_{A_\alpha}  \right\|_{L_\infty}
\hsp{3}-\hsp{10} 
\left\|\sum_{n=1}^m a_n \left( \at_n - f_n\right)
     1_{A_\alpha}   \right\|_{L_\infty}
\\ &\ge\hsp{5}
\sum_{n=1}^m |a_n|  \hsp{3}-\hsp{3} \sum_{n=1}^m \e_n |a_n| 
\hsp{5}=\hsp{5}
\sum_{n=1}^m (1-\e_n) |a_n| \ . 
\end{align*}
So  $(\e_n)$ and $(x_n)$ do   indeed satisfy (A1).

To show that (a) implies (f), let (a) hold. 
Then we have the situation depicted in~\eqref{eq: cd}.   
For $t \in [0,1]$, let $\delta_t\in C^*$ denote the
point mass measure at~$t$. 
Since  $T^*$ is      
a  weak-star continuous isometric embedding
\begin{equation*} 
M \defeq  \{T^*(\delta_t) \colon t\in [0,1]\}
\hsp{4} \subset \hsp{4}  B(Y^*) \ , 
\end{equation*}
equipped with the weak-star topology of $Y^*$,
is homeomorphic to $[0,1]$ and 
is isometrically equivalent to
the usual basis of $\ell_1([0,1])$.
By the Hahn-Banach Theorem
$
j^*(B(X^*))=B(Y^*)   
$  
and so 
\begin{equation*} 
N \defeq j^{*-1}(M)\cap B(X^*)
\end{equation*} 
is weak-star compact and satisfies $j^*(N)=M$.
By Fact~\ref{topology} there exists 
\begin{equation*}
\Gamma=\{ n_t \colon t \in [0,1]\} \subset N
\end{equation*}
such that  
$
j^*(n_t) = T^*(\delta_t)
$  
and such that $\Gamma$ is dense-in-itself
in the weak-star topology on $X^*$.  
Since  for any finite set $\{ a_t \}_{t\in F}$ 
of scalars 
\begin{alignat*}{2}
\sum_{t\in F} \left| a_t \right|
\hsp{3}&=\hsp{3} 
\left\| \sum_{t\in F} a_t T^* \delta_t \right\|_{Y^*} 
&\hsp{3}&=\hsp{3} 
\left\| \sum_{t\in F} a_t j^*  n_t \right\|_{Y^*} 
\\ \hsp{3}&\le\hsp{3} 
\left\| \sum_{t\in F} a_t n_t \right\|_{X^*}
&\hsp{3}&\le\hsp{3}
\sum_{t\in F} \left| a_t \right| \ , 
\end{alignat*}
the set $\Gamma$ is
isometrically
equivalent to the usual basis of $\ell_1([0,1])$.

To see that (f) implies (e),  
let $X$ satisfy (f). Then by Lemma~\ref{lem: sep}, 
there is a separable subspace $X_0$  of $X$ which   
satisfies (f). From   Fact~\ref{thm: hpr}       
and the fact that $C^*(\Delta)$ is linearly isometric to~$C^*$, 
it follows that~$X_0$ satisfies (e).  
The equivalence of (a) and (e) gives  
that~$X$  also satisfies~(e).

Thus  (a) through (f) are equivalent.  
The fact that 
$C^*(\Delta)$ is linearly isometric to~$C^*$ gives that 
(g) implies (e).  
That (f) implies (g) when $X$ is separable 
is due to Rosenthal:  Fact~\ref{thm: hpr}.     
\end{proof}

\begin{remark}
Without the added assumption of separability, 
(g) is not equivalent to the other conditions. 
Clearly, $\ell_\infty$ satisfies conditions~(a) 
through~(f).   
But by a result of Grothendieck \cite{G}, a separable 
quotient of $\ell_\infty$ is reflexive and so~$\ell_\infty$ 
does not satisfy~(g). 
\end{remark}

\begin{remark}  
A complemented isomorphic version of Theorem~\ref{th: mainresult} 
is due to Hagler and Stegall \cite[Theorem~1]{HS}. 
A $K$-complemented isometric version of Theorem~\ref{th: mainresult} 
is due to Hagler (\cite[Theorem~2.13]{H1} or \cite{H3}).   
\end{remark}

\begin{remark} Many  Banach
spaces (and their subspaces)
which arise naturally in analysis contain an abundance of
asymptotically isometric copies of~$\ell_1$: for example,
Carothers, Dilworth and Lennard~\cite{CDL} proved that
every nonreflexive subspace of the Lorentz space $L_{w,1}(0,\infty)$
contains asymptotically isometric copies of $\ell_1$ whenever
the weight~$w$ satisfies very mild regularity conditions.
On the other hand, $L_{w,1}$ 
does {\it not} contain an isometric copy
of the~$2$-dimensional space $\ell_1^2$ whenever 
$w$ is strictly decreasing \cite{CDT}. \end{remark}

\begin{remark} Theorem~\ref{th: mainresult} improves a recent result of
Shutao Chen and Bor-Luh Lin \cite{CL} who proved that $X$ contains
an asymptotically isometric copy of~$\ell_1$ whenever $X^*$
contains an isometric copy of $\ell_\infty$.
\end{remark}

Dowling, Johnson, Lennard and Turett~\cite{DJLT} 
gave some concrete examples of
equivalent norms on $\ell_1$
such that the corresponding    Banach spaces do  not  
 contain asymptotically isometric copies of $\ell_1$.
Theorem~\ref{th: mainresult} easily yields other
equivalent norms on $\ell_1$ with this property.

\begin{corollary} 
Let $(\gn)$   
be a sequence of non-zero  scalars with $\|(\gn)\|_2<\gnn$. 
Then   the  Banach space~$(\ell_1,\|\cdot\|_1')$, where
\begin{equation*}
\|(a_n)\|_1'  \defeq   
\inf\left\{\left[\,\left\| \left(a_n + b_n \right)\right\|_1^2 ~+~
\left\| \left( \gn^{-1} b_n \right) \right\|_2^2 \,\right]^{\frac{1}{2}} 
\colon \left( \gn^{-1} b_n \right) \in \ell_2 \right\} \ , 
\end{equation*} 
 does not contain asymptotically isometric copies of $\ell_1$  
and 
\begin{equation}
\label{eq:equiv}
\left(1+\e^2\right)^{-\frac{1}{2}} 
\left\| \left( a_n \right) \right\|_1  
~\leq~
\left\| \left( a_n \right) \right\|_1'
~\leq~
\left\| \left( a_n \right) \right\|_1  
\end{equation} 
for each $\left( a_n \right) \in \ell_1$. 
\end{corollary}  
\begin{proof} 
We exhibit $(\ell _{1},\Vert \cdot \Vert _{1}^{\prime })$ as a quotient
space of $X=\ell _{1}\oplus _{2}\ell _{2}$ with its usual norm
\begin{equation*}
\Vert ((a_{n}),(b_{n}))\Vert_{X}  
=\left[\, \left\| \left( a_n \right)\right\|_{1}^2 ~+~
\left\| \left( b_n \right)\right\|_{2}^2 \,\right]^{1/2}.
\end{equation*}
Let 
\begin{alignat*}{2}
Y &\defeq &&\{((a_{n}),(b_{n}))\in X ~\colon~ b_{n}= -\gn^{-1}a_{n}\} \\
&\hsp{5}\overset{\text{\rm note}}{=}\hsp{5}
&&\{((a_{n}),(-\gn^{-1}a_{n} ))  ~\colon~ (\gn^{-1}a_{n})\in\ell_2\} \ . 
\end{alignat*} 
Since  each element of $X/Y$
has a representative of the 
form $\left( \left( a_{n}\right) ,0\right)$,   
\begin{multline*}
\Vert ((a_{n}),0)+Y\Vert _{X/Y} \hsp{5} = \\
\inf \left\{\, 
\left\| \left( \left(a_n\right), 0 \right) ~+~ 
\left(\left( b_n\right) , \left(-\gn^{-1}b_n\right)\right) \right\|_X 
~\colon~ \left(\gn^{-1}b_n\right)\in\ell_2 \,\right\}
\\ = \hsp{5} \Vert (a_{n})\Vert _{1}^{^{\prime }} \ . 
\end{multline*}
Thus  $X/Y$ is isometrically isomorphic to 
$(\ell _{1},\Vert \cdot \Vert_{1}^{\prime })$.   
 
Observe that $X^{\ast }=\ell _{\infty }\oplus _{2}\ell _{2}$ with its usual
norm 
\begin{equation*}
\Vert ((c_{n}),(d_{n}))\Vert_{X^*}  
=\left[\, \left\| \left( c_n \right)\right\|_{\infty}^2 ~+~
\left\| \left( d_n \right)\right\|_{2}^2 \,\right]^{1/2}   
\end{equation*}
and 
$$
Y^{\perp} 
\hsp{2}=\hsp{2} 
\{((c_{n}),(d_{n}))\in X^* ~\colon~ d_{n}=  \gn c_{n}\} 
\hsp{2}=\hsp{2}  
\{((c_{n}),( \gn c_{n} ))  ~\colon~ ( c_n)\in\ell_\infty\}~. 
$$      
It follows 
that $Y^{\perp }$ is isometrically isomorphic to $(\ell _{\infty },\Vert
\cdot \Vert _{\infty }^{\prime })$ where 
\begin{equation}
\label{eq:equivc}
\Vert (c_{n})\Vert_{\infty }^{\prime }
=\left[\, \Vert (c_{n})\Vert _{\infty}^{2}
~+~ 
\left\| \left( \gn c_n \right) \right\|_2^2 \,\right]^{1/2}   
\end{equation}  
and the mapping 
\begin{equation*} 
i \colon \left(\ell_\infty  ,\Vert \cdot \Vert _{\infty }^{\prime }\right) 
\to 
\left(\ell_1  ,\Vert \cdot \Vert _{1}^{\prime }\right)^*
\end{equation*}
given by $\left(i\left(c_n\right)\right)\left(a_n\right) = \sum_n a_n c_n$ 
is an isometry.  
Since the norm $\Vert \cdot \Vert _{\infty }^{\prime }$ 
is strictly convex, 
the space $(\ell _{\infty },\Vert \cdot \Vert _{\infty }^{\prime })$ 
does not contain an isometric copy of~$L_{1}$. 
Thus by Theorem~\ref{th: mainresult}  
the space  $(\ell _{1},\Vert \cdot \Vert _{1}^{\prime })$ does not contain
asymptotically isometric copies of~$\ell _{1}$.  From \eqref{eq:equivc} 
it follows that 
\begin{equation*} 
\left\| \left( c_n \right)\right\|_\infty 
~\le~
\left\| \left( c_n \right)\right\|_\infty^\prime 
~\le~
\sqrt{1+\e^2}\, \left\| \left( c_n \right)\right\|_\infty 
\end{equation*} 
for each $\left( c_n \right)\in\ell_\infty$ 
and so  \eqref{eq:equiv} holds  by duality. 
\end{proof}  
Finally,  Alspach's~\cite{A} example of
 an isometry $T \colon K\rightarrow K$ without a fixed point, where $K$
is a certain weakly compact convex subset of $L_1$,
yields an obvious corollary.
\begin{corollary}     If  $X$ contains asymptotically isometric
copies of~$\ell_1$, then there exists a (nonempty)
 weakly compact convex subset
 $K$ of $ X^*$ and an isometry $T\colon K \rightarrow K$ without a fixed point.
\end{corollary}
\begin{Acknowledgement} 
The authors thank  
Patrick Dowling, William B.\ Johnson, and Haskell P.\ Rosenthal 
for  their helpful comments. 
H.~P.~Rosenthal 
suggested the term p.i.\ $\ell_1$-sequence 
and  his helpful conversations~\cite{HPR}   
led us to the proof of Lemma~\ref{lem:pi}. 
\end{Acknowledgement}

\end{document}